\documentclass[10pt]{article}
\usepackage{amssymb,amsmath,amsthm,latexsym}
\usepackage{amsfonts}
\usepackage[twlrm]{rawfonts}
\usepackage[english]{babel}
\usepackage[latin1]{inputenc}
\usepackage{fancyhdr}
\usepackage{indentfirst}

\newtheorem{teo}{Theorem}

\newtheorem{lem}[teo]{Lemma}

\newenvironment{dem}[1][Proof]{\noindent\textbf{#1.} }{\hfill \rule{0.5em}{0.5em}}

\newcommand{\N}{\mathbb{N}}
\newcommand{\Z}{\mathbb{Z}}
\newcommand{\R}{\mathbb{R}}

\newcommand{\ddd}{\displaystyle}

\newtheorem{ob}{Remark}[section]

\newtheorem{Af}{Claim}[section]

\begin{document}
	
\setlength{\baselineskip}{6.5mm} \setlength{\oddsidemargin}{8mm}
\setlength{\topmargin}{-3mm}

\title{\bf Choquard Equations with Mixed Potential}

\author{Marco A. S. Souto\thanks{M. A. S. Souto was partially supported by CNPq/Brazil
			305384/2014-7 and INCT-MAT, marco@dme.ufcg.edu.br}\, \, and \, \  Romildo N. de Lima \thanks{ R. N. de Lima was partially supported by CAPES/Brazil, romildo@dme.ufcg.edu.br}\,\,\,\,\,\,\vspace{2mm}
		\and {\small  Universidade Federal de Campina Grande} \\ {\small Unidade Acadêmica de Matemática} \\ {\small CEP: 58429-900, Campina Grande - PB, Brazil}\\}

\date{}

\maketitle

\begin{abstract}
	In this paper, we study the following class of nonlinear Choquard equation,
	$$-\Delta u+a(z)u=K(u)f(u)\quad \text{in}\quad \R^N,$$
where $\R^N=\R^L\times\R^M$, $L\geq2$, $K(u)=|.|^{-\gamma}*F(u)$, $\gamma\in(0,N)$, $a$ is a continuous real function and $F$ is the primitive function of $f$. Under some suitable assumptions mixed on the potential $a$. We prove existence of a nontrivial solution for the above equation. \vspace{0.5cm}


\end{abstract}
	
\section{Introduction and main results}

The aim of this paper is to study the existence of nontrivial solutions for the following nonlinear Choquard equation family

\[
(P)\left\{
\begin{array}{lcl}
-\Delta u+a(z)u =K(u)f(u),\quad \R^{N} \\
u\in H^{1}(\R^{N}),\quad u>0\quad\text{em $\R^{N}$} 
\end{array}
\right.
\] 
where $K(u)=|.|^{-\gamma}*F(u)$, $\gamma\in(0,N)$, $a$ is a nonnegative continuous real function and $F$ is the primitive function of $f$.

This problem comes from observation of waves solutions for a nonlinear Schrodinger equation of the kind

\begin{equation}
i\partial_{t}\Psi=-\Delta\Psi+W(z)\Psi-(Q(z)*|\Psi|^{q})|\Psi|^{q-2}\Psi, \quad\text{in}\quad \R^N.
\end{equation}

In this context, $W$ is the external potential and $Q$ is the response function possesses information on the mutual interaction between the bosons. This type of nonlocal equation is known to influence the propagation of electromagnetic waves in plasmas \cite{Berge-Couairon} and also plays an important role in the theory of Bose-Einstein condensation \cite{Dalfovo-Giorgini}. It is clear that $\Psi(z,t)=u(z)e^{-iEt}$ solves the evolution equation $(1)$ if, and only if, $u$ solves

\begin{equation}
-\Delta u+a(z)u=(Q(z)*|u|^{q})|u|^{q-2}\quad\text{in}\quad \R^N,
\end{equation}
with $a(z)=W(z)-E$.

When the response function is the Dirac function, i.e. $Q(z)=\delta(z)$, the nonlinear response is local and we have the Schrödinger equation

\begin{equation}
-\Delta u+a(z)u=|u|^{q-2}u\quad\text{in}\quad \R^N.
\end{equation}
This equation has been studied extensively under various hypotheses on the potentials and the nonlinearities. We may refer to \cite{Alves-Souto-Montenegro}, \cite{Alves-Souto}, \cite{Bartsch-Pankov-Wang}, \cite{Byeon-Wang}, \cite{Berestycki-Lions} and the references therein.

In this paper, we study the existence of nontrivial solutions for a class of Schrödinger equation with nonlocal type nonlinearities, that is, the response function $Q$ in $(P)$ is of Coulomb type, for example $|z|^{-\gamma}$, then we arrive at the Choquard-Pekar equation,

\begin{equation}
-\Delta u+a(z)u=\left(\frac{1}{|z|^{\gamma}}*|u|^{q}\right)|u|^{q-2}\quad\text{in}\quad \R^N.
\end{equation}

We know that, most of the existing papers consider the existence and property of the solutions for the nonlinear Choquard equation $(P)$ with constant potential. For example: in \cite{Lieb}, Lieb proved the existence and uniqueness, up to translations, of the ground state solution to equation $(4)$. Later, in \cite{Lions}, Lions showed the existence of a sequence of radially symmetric solutions to this equation. Involving the properties of the ground state solutions, Ma and Zhao \cite{Ma-Zhao} considered the generalized Choquard equation $(4)$ for $q\geq2$, and they proved that every positive solution of $(4)$ is radially symmetric and monotone decreasing about some point.

Involving the problem with nonconstant potentials, we have

\begin{equation}
-\Delta u+a(z)u=\left(\frac{1}{|z|^{\gamma}}*F(u)\right)f(u)\quad\text{in}\quad \R^N,
\end{equation}
where $a$ is a continuous periodic function with $\inf_{\R^N}a(z)>0$, noticing that the nonlocal term is invariant under translation, it is possible to prove an existence result easily by applying the Mountain Pass Theorem, see \cite{Ackermann} for example. In \cite{Alves-Figueiredo-Yang}, Alves, Figueiredo and Yang, they made a very interesting work for the case generalized Choquard equation with vanishing potential.

Looking at the various works cited above and others, we observe the lack of existence of results, for mixed potentials with different characteristics in each entry, i.e., considering $\R^N=\R^L\times\R^M$ and $z\in\R^N$ as $z=(x,y)$, $x\in\R^L$ and $y\in\R^M$, $a(x,y)$ has different characteristics for each variable. This led us to seek solution to some kinds of mixed potential.

In all cases we study, $f:\R\longrightarrow\R$ is continuous and satisfies:

$(f_{0})$ $|f(s)|\leq C_{0}(|s|^{q_{1}-1}+|s|^{q_{2}-1})$, where $q_{1},q_{2}>1$ with $\frac{2N-\gamma}{N}<q_{1}\leq q_{2}<\frac{2N-\gamma}{N-2}$;

$(f_{1})$ $\frac{f(t)}{t}$ is increasing and unbounded in $t>0$;

$(f_{2})$ $f(t)>0$ in $t>0$ and $f(t)=0$ in $t\leq0$.

The first mixed potential is the symmetric-coercive  type, i.e., $a:\R^L\times\R^M\longrightarrow\R$ continuous, where $L\geq2$, such that,

$(a_{0})$ There exists $R>0$ such that $a(x,y)\geq a_{0},\quad\text{for all}\quad z\in B_{R}^{c}(0),\quad z=(x,y)$ and $a(z)\geq0,\quad\text{for all}\quad z\in\R^{N}$,

$(a_{1})$ $a(x,y)\longrightarrow\infty$, when $|y|\longrightarrow\infty$ uniformly for $x\in\R^L$;

$(a_{2})$ $a(x,y)=a(x',y)$ for all $x,x'\in\R^L$ with $|x|=|x'|$ and all $y\in\R^M$.

In this case, the obtained solution $u\in H^{1}(\R^N)$ is such that, $u(x,y)=u(x',y)$, always that $x,x'\in\R^L$ with $|x|=|x'|$ and all $y\in\R^M$. To do this, we will make use of \cite{Cordeiro-Souto-O}, wherein the authors prove a compactness embedding lemma and a principle of symmetric criticality. In fact, we prove:

\begin{teo}
	Assume $(a_{0})-(a_{2})$ and $(f_{0})-(f_{2})$. Then, problem $(P)$ has a positive solution.
\end{teo}

The second mixed potential is the periodic-coercive type, i.e., $a:\R^L\times\R^M\longrightarrow\R$ continuous, such that, satisfies $(a_0)$, $(a_1)$ and

$(a_{3})$ $a(x,y)=a(x+p,y)$, for all $x\in\R^L,y\in\R^M$ and $p\in\Z^L$.

In fact, we prove:

\begin{teo}
	Assume $(a_{0}),(a_{1})$ and $(a_{3})$ for $a$ and $(f_{0})-(f_{2})$ for $f$. We have that, the problem $(P)$ has a positive solution in the  level of mountain pass.
\end{teo}

The third and last mixed potential we work is the asymptotically periodic-coercive type, i.e., $a:\R^L\times\R^M\longrightarrow\R$ continuous, such that, satisfies $(a_0)$, $(a_1)$ and

$(a_4)$ There exists a potential $a_p:\R^L\times\R^M\longrightarrow\R$ continues such that, satisfies $(a_0)$, $(a_1)$, $(a_3)$ and

\begin{enumerate}
	\item There exists $\mu>0$ such that,
	$$a(z)\lneq a_p(z)\leq\mu a(z),\quad\text{for all}\quad z\in\R^N$$
	where, "$a(z)\lneq a_p(z)$", means that, there exists $\Omega\subset\R^N$ with $|\Omega|>0$, such that, $a(z)<a_p(z)$ for all $z\in\Omega$;
	\item $|a(x,y)-a_p(x,y)|\longrightarrow0$, when $|y|\longrightarrow\infty$, uniformily for $x\in\R^L$.
\end{enumerate}

In fact, we prove:

\begin{teo}
	Assume $(a_{0}),(a_{1})$ and $(a_{4})$ for $a$ and $(f_{0})-(f_{2})$ for $f$. Then, problem $(P)$ has a positive solution.
\end{teo}

\begin{ob}
	In our study, it was not necessary that the nonlinearity $f$ verifies the Ambrosetti-Rabinowitz type superlinear condition for nonlocal problem, see \cite{Ambrosetti-Rabinowitz}, that is, there exists $\theta>2$, such that
	$$0<\theta F(s)\leq2f(s)s,\quad\text{for all}\quad s>0,$$
	we often found this hypothesis in work on this subject. Under this condition, all results \textbf{theorem 1}, \textbf{theorem 2} and \textbf{theorem 3} are true replacing $(f_1)$ by Ambrosetti-Rabinowitz condition.
\end{ob}

\begin{ob}
	A common hypothesis in the work when do not have the Ambrosetti-Rabinowitz condition, is the condition
	$$\lim_{t\longrightarrow\infty}\frac{F(t)}{t^2}=\infty,$$
	but, it is easy to see that $(f_1)$ implies the above condition.
\end{ob}

\textbf{Notations}

We fix the following notations, which will use from now on.

\begin{itemize}
	
	\item $z\in\R^N=\R^L\times\R^M$ is given by $z=(x,y)$, such that, $x\in\R^L$ and $y\in\R^M$.
	
	\item $B_{R}(z)$ denotes the ball centered at the $z$ with radius $R>0$ in $\R^N$.
	
	\item $L^{s}(\R^N)$, for $1\leq s\leq\infty$, denotes the Lebesgue space with the norms
	$$|u|_{s}=\left(\int_{\R^N}|u|^{s}dz\right)^{\frac{1}{s}}$$
	or
	$$|u|_{s}=\left(\int_{\R^L\times\R^M}|u|^{s}dxdy\right)^{\frac{1}{s}}.$$
	
	\item $C_{0}^{\infty}(\R^N)$ denotes the space of the functions infinitely differentiable with compact support in $\R^N$.
	
	\item We denote the inner product of $H^{1}(\R^N)$ by
	$$(u,v)_{H^1}=\int_{\R^N}\nabla u\nabla v+uvdz$$
	and the norm
	$$\|u\|_{H^{1}}=\left(\int_{\R^N}|\nabla u|^{2}+|u|^{2}dz\right)^{\frac{1}{2}}.$$
	
	\item From the assumptions on $a$, it follows that the subspace
	$$E_{a}=\{u\in H^{1}(\R^N);\quad \int_{\R^N}a(z)|u|^{2}dz<\infty\}$$
	is a Hilbert space with norm defined by
	$$\|u\|:=\left(\int_{\R^N}|\nabla u|^{2}+a(z)|u|^{2}dz\right)^{\frac{1}{2}}$$
	and $E_{a}\hookrightarrow H^{1}(\R^N)$, continuously.
\end{itemize}
	
\section{Common properties of the problem with Mixed Potential}

We would like to make some comments on the assumptions involving the nonlinearity $f$. We intend to use variational methods, this way, we must have:
$$\left|\int_{\R^{N}}(|z|^{-\gamma}*F(u))F(u)dz\right|<\infty,\quad\text{for all}\quad u\in E_{a}.$$

To see that above property occurs, it is very important to recall the Hardy-Littlewood-Sobolev inequality, found in \cite{Lieb-Loss}, which will be frequently used in the paper.

\begin{teo}
	Let $p,r>0$ and $\gamma\in(0,N)$ with $1/p+\gamma/N+1/r=2$. If $f\in L^{p}(\R^{N})$ and $h\in L^{r}(\R^{N})$, then there exists a sharp constant $C:=C(p,N,\gamma,r)>0$, independent of $f$ and $h$, such that
	$$\ddd\int_{\R^{N}}\ddd\int_{\R^{N}}\frac{f(x)h(y)}{|x-y|^{\gamma}}dydx\leq C|f|_{p}|h|_{r}.$$
\end{teo}

By $(f_{0})$,
$$|F(u)|\leq C_{0}(|u|^{q_{1}}+|u|^{q_{2}})\Longrightarrow |F(u)|^{\frac{2N}{2N-\gamma}}\leq C_{0}(|u|^{q_{1}\frac{2N}{2N-\gamma}}+|u|^{q_{2}\frac{2N}{2N-\gamma}})$$
as,
$$2<q_{1}\frac{2N}{2N-\gamma}\leq q_{2}\frac{2N}{2N-\gamma}<2^{*}$$
thence, $F(u)\in L^{\frac{2N}{2N-\gamma}}(\R^{N})$. Now, note that,
$$\frac{2N-\gamma}{N}+\frac{\gamma}{N}=2,$$
thereby, by Hardy-Littlewood-Sobolev inequality,
$$\left|\int_{\R^{N}}K(u)F(u)dz\right|\leq C_{0}|F(u)|_{\frac{2N}{2N-\gamma}}^{2}<\infty.$$

From the above commentaries, the Euler-Lagrange functional $I:E_{a}\longrightarrow\R$ associated to $(P)$ given by
$$I(u):=\frac{1}{2}\int_{\R^{n}}|\nabla u|^{2}+a(z)|u|^{2}dz-\frac{1}{2}\int_{\R^{n}}K(u)F(u)dz,$$
is well defined and belongs to $C^1$ with its derivative given by
$$I'(u)\varphi=\int_{\R^N}\nabla u\nabla\varphi+a(z)u\varphi dz-\int_{\R^N}K(u)f(u)\varphi dz,$$
for all $u,\varphi\in E_{a}$. Thus, it is easy to see that all the solutions of $(P)$ correspond to critical points of the energy functional $I$.

We have that, $I$ verifies the mountain pass geometry, through of arguments well know in the literature.

\begin{lem}
	Assume $\gamma\in(0,N)$, $(a_{0})$ and $(f_{0})-(f_{2})$. Then,	
	
	$(1)$ There exist $\rho,\delta_{0}>0$ such that $I|_{S_{\rho}}\geq \delta_{0}>0$;
	
	$(2)$ There exist $\varphi\in E_{a}$ with $\|\varphi\|>\rho$ such that $I(\varphi)<0$.
	
\end{lem}

By the mountain pass theorem, see \cite{Ekeland}, there is a Cerami sequence $(u_{n})\subset E_{a}$, such that
$$I(u_{n})\longrightarrow c_{a}\quad\text{e}\quad I'(u_{n})(u_{n})\longrightarrow0,$$
where
$$c_{a}:=\inf_{\alpha\in\Gamma}\max_{t\in[0,1]}I(\alpha(t))>0$$	
and
$$\Gamma:=\{\alpha\in C([0,1],E_{a});\quad \alpha(0)=0\quad\text{e}\quad I(\alpha(1))<0\}.$$

The following lemma will be useful to prove our results.

\begin{lem}
	Let $(u_{n})\subset E_{a}$ such that $(I(u_{n}))$ is bounded and $\|u_{n}\|\longrightarrow\infty$, then, $w_{n}=\frac{u_{n}}{\|u_{n}\|}$ is such that $w_{n}\rightharpoonup w$ in $E_{a}$, where $w\leq0$ a.e. in $\R^{N}$.
\end{lem}

\begin{dem}
	Since $(w_{n})$ is bounded in $E$, there are $w\in E_{a}$ and a subsequence, still denote by $(w_{n})$, such that $(w_{n})\rightharpoonup w$ in $E_{a}$. For all $R>0$, we define
	$$\Omega:=\{z\in\R^{N};\quad |z|\leq R\quad\text{e}\quad w(z)>0\}.$$
	We claim that $|\Omega|=0$, for all $R>0$. Suppose by contradiction that $|\Omega|>0$, for some $R>0$. Thus,
	$$I(u_{n})=c+o_{n}(1)\Longrightarrow \int_{\R^{N}}K(u_{n})\frac{F(u_{n})}{\|u_{n}\|^{2}}dz=1+o_{n}(1)$$	
	from where it follows that,
	$$\int_{\Omega}K(u_{n})\frac{F(u_{n})}{\|u_{n}\|^{2}}dz\leq 1+o_{n}(1).$$
	By \textbf{Remark 1.2}, for all $M>0$, there exists $\delta>0$ such that
	$$\frac{F(s)}{s^{2}}\geq M,\quad\text{for all}\quad s\geq\delta$$
	thence, for
	$$G_{n}:=\{z\in\R^{N};\quad u_{n}(z)\geq\delta\}$$
	we have to,
	$$\int_{\Omega\cap G_{n}}K(u_{n})\frac{F(u_{n})}{u_{n}^{2}}w_{n}^{2}dz\leq 1+o_{n}(1)$$
	where,
	$$M\int_{\Omega\cap G_{n}}K(u_{n})w_{n}^{2}dz\leq 1+o_{n}(1).$$
	
	Note that, for $z\in\Omega$, from $n$ large enough $u_{n}(z)\geq\delta$, because, $w_{n}(z)\longrightarrow w(z)$ and $u_{n}(z)=\|u_{n}\|w_{n}(z)$. Now, for $z\in\Omega$, $n$ large enough,
	$$K(u_{n})(z)\geq F(\delta)\int_{\R^{N}}\chi_{\Omega\cap G_{n}}(z')|z-z'|^{-\gamma}dz'$$
	and so
	$$\liminf K(u_{n})(z)\geq F(\delta)>0.$$
	Therefore,
	$$M\liminf \int_{\Omega\cap G_{n}}K(u_{n})w_{n}^{2}dz\leq1$$
	by Fatou lemma,
	$$MF(\delta)\int_{\Omega}w^2dz\leq1,\quad\text{for all}\quad M>0,$$
	consequently, $|\Omega|=0$. 
\end{dem}

\section{Symmetric-Coercive Case}

In this case, due to the lack of compactness, we restrict $I$ to a subspace of $E_{a}$, given by
$$E=:\{u\in E_{a};\quad u(x,y)=u(x',y),|x|=|x'|\}$$
thus, $E$ is Hilbert space under the scalar product,
$$(u,v)_{E}:=\int_{\R^{N}}(\nabla u\nabla v+a(z)uv) dz.$$

Thus, we have the lemma:

\begin{lem} $E$ is continuously immersed in $L^{s}(\R^N)$ if $s\in[2,2^*]$ and compactly immersed if $s\in(2,2^*)$.

Indeed, first of all, we are going to prove that if condition $(a_0)$ is valid then the Banach space $E$ is continuously immersed in $L^{s}(\R^N)$ for all $s\in[2,2^*]$. Notice that $(a_0)$ yields that,
$$\int_{\R^N}|u|^2dz=\int_{|z|\leq R}|u|^2dz+\int_{|z|>R}|u|^2dz\leq\int_{|z|\leq R}|u|^2dz+a_0^{-1}\int_{|z|>R}a(z)|u|^2dz\leq C\|u\|_{E}^{2},$$
where we use the continuity of the Sobolev embedding for
bounded domains. On the other hand, the Sobolev-Gagliardo-Nirenberg
inequality asserts that there exists positive constant $S$ such that
$$\int_{\R^N}|u|^{2^*}dz\leq S\int_{\R^N}|\nabla u|^{2}dz.$$
Therefore, from inequalities above, we have the continuity of
the embedding for $s=2$ and $s=2^*$. The continuity of the immersion for a fixed $s\in(2,2^*)$, follows from the following interpolation
inequality
$$|u|_s\leq|u|_{2}^{1-t}|u|_{2^*}^{t}$$
where $t\in[0,1]$ is such that $1/s=(1-t)/2+t/2^*$.

Now, assuming $(a_0)-(a_2)$, we are going to prove the compactness
of the embedding of the spaces $E$ in $L^s(\R^N)$ for $s\in(2,2^*)$. Let $(u_n)\subset E$ bounded, so, $(u_n)$ is bounded in $H^1(\R^N)$, thus, $u_n\rightharpoonup u$ in $H^1(\R^N)$. Setting $v_n=u_n-u$, we have, by Lion's Lemma:

$(A)$ There exist $(z_n)\subset\R^N$ and $\rho,R>0$, such that
$$\int_{B_R(z_n)}v_n^2dz\geq\rho>0\quad\text{for all}\quad n\in\N,$$
or

$(B)$ $v_n\longrightarrow0$ in $L^q(\R^N)$, for all $q\in(2,2^*)$.

If $(A)$ occurs, we have that $(y_n)$ is bounded in $\R^M$. Indeed $(v_n)$ is bounded, thus, there exists $\Lambda>0$ such that,
$$\frac{\|v_n\|^2}{\Lambda}<\rho,\quad\mbox{for all}\quad n\in\N$$
so, noticing that, there exists $\lambda>0$ such that,
$$a(x,y)\geq\Lambda,\quad\mbox{for all}\quad y\in B_\lambda^c(0)$$
if $(y_n)$ is \textbf{unbounded}, for $n$ large enough, then
$$(x,y)\in B_R(x_n,y_n)\Longrightarrow|y|\geq\lambda\Longrightarrow a(x,y)\geq\Lambda$$
thus,
$$0<\rho\leq\int_{B_R(x_n,y_n)}|v_n|^2dz\leq\frac{1}{\Lambda}\int_{B_R(x_n,y_n)}a(z)|v_n|^2dz\leq\frac{\|v_n\|^2}{\Lambda}<\rho,$$
which is a contradiction, therefore, $(y_n)$ is bounded. And so, there exists $\lambda>0$ such that $|y_n|\leq\lambda$ for all $n\in\N.$ Therefore, for $\overline{R}=R+\lambda$, we have
$$B_{R}(x_n,y_n)\subset B_{\overline{R}}(x_n,0)\Rightarrow \int_{B_{\overline{R}}(x_n,0)}|v_n|^2dz\geq\rho>0,\quad\mbox{for all}\quad n\in\N.$$

If $(x_n)$ is unbounded, without loss of generality we can assume that, for every $n\in\N$, $|x_n|\geq n\overline{R}$. Thus, there are at least $3n$ balls of ray $\overline{R}$, disjointed and centered at points $(\overline{x},0)$, for $|\overline{x}|=|x_n|$. We denote by J such finite set of centers, and so
$$\int_{B_{\overline{R}}(x_n,0)}|v_{n}|^{2}dz=\int_{B_{\overline{R}}(\overline{x},0)}|v_{n}|^{2}dz,\quad\mbox{for all}\quad (\overline{x},0)\in J$$
thus,
$$\int_{\R^N}|v_n|^2dz\geq\sum_{(\overline{x},0)\in J}\int_{B_{\overline{R}}(\overline{x},0)}|v_{n}|^{2}dz\geq3n\rho$$
and this is a contradiction because, $(v_n)$ is bounded. Therefore, $(x_n)$ is bounded. Which also creates an absurd. And so, $(B)$ is valid and we have the compact immersion.

\end{lem}

The proof of these immersions can be made using \textbf{Lemma III.2} pp. 321 of \cite{Lions2}, but, we choose this, because it is simpler.

We seek critical point of $I|_{E}$, and by principle of symmetric criticality in \cite{Cordeiro-Souto-O}, this point is critical in $I:E_a\longrightarrow\R$.

It is very important to note that, the lemmas above are valid, replacing $E_{a}$ by $E$. Now,

\begin{lem}
	If $v_{n}\rightharpoonup v$ in $E$, then
	$$\lim_{n\longrightarrow\infty}\int_{\R^N}K(v_{n})F(v_{n})dz=\int_{\R^N}K(v)F(v)dz$$
	and
	$$\lim_{n\longrightarrow\infty}\int_{\R^N}K(v_{n})f(v_{n})v_{n}dz=\int_{\R^N}K(v)f(v)vdz$$
\end{lem}

\begin{dem}
	Note that, $v_{n}\rightharpoonup v$ in $E$. Thus, from compact immersion of $E$ in $L^s(\R^N)$ for $s\in(2,2^*)$, $v_{n}(z)\longrightarrow v(z)$ and $F(v_{n}(z))\longrightarrow F(v(z))$, $a.e.$ in $\R^N$ and $v_{n}\longrightarrow v$ in $L^t(\R^N)$, $t\in(2,2^*)$.
	Now, recalling that,
	$$|F(v_{n})|^{\frac{2N}{2N-\gamma}}\leq C_{0}(|v_{n}|^{q_{1}\frac{2N}{2N-\gamma}}+|v_{n}|^{q_{2}\frac{2N}{2N-\gamma}})$$
	and
	$$2<q_{1}\frac{2N}{2N-\gamma}\leq q_{2}\frac{2N}{2N-\gamma}<2^{*},$$
	we have, using the dominated convergence theorem and the fact that $\Sigma:L^{\frac{2N}{2N-\gamma}}(\R^N)\longrightarrow L^{\frac{2N}{\gamma}}(\R^N)$, given by $\Sigma(w):=|.|^{-\gamma}*w$ is a linear bounded operator, that
	$$K(v_{n})\longrightarrow K(v),\quad\text{em}\quad L^{\frac{2N}{\gamma}}(\R^N).$$
	Using similar arguments, we show that
	$$f(v_{n})v_{n}\longrightarrow f(v)v\quad\text{and}\quad f(v_{n})v\rightharpoonup f(v)v\quad\text{in}\quad L^{\frac{2N}{2N-\gamma}}(\R^N).$$
	Therefore,
	$$\int_{\R^N}K(v_{n})f(v_{n})v_{n}dz\longrightarrow\int_{\R^N}K(v)f(v)vdz$$
	and
	$$\int_{\R^N}K(v_{n})f(v_{n})vdz\longrightarrow\int_{\R^N}K(v)f(v)vdz.$$	
\end{dem}

\begin{lem}
	Let $(u_{n})\subset E$ such that $I(u_{n})\longrightarrow c_{a}$ and $I'(u_{n})(u_{n})\longrightarrow0$, then, $(u_{n})$ is bounded in $E$.
\end{lem}

\begin{dem}
	Indeed, otherwise, $\|u_{n}\|\longrightarrow\infty$. Setting, $w_{n}=u_{n}/\|u_{n}\|$, we have from \textbf{Lemma 6}
	$$w_{n}\rightharpoonup w\quad\text{in}\quad E,$$
	with $w\leq0$. We can check that, \textbf{Lemma 8} implies that
	$$\lim_{n\longrightarrow\infty}\int_{\R^N}K(w_{n})f(w_{n})w_{n}dz=0.$$
	Now, note that, for all $n\in\N$,
	$$\text{there exists}\quad t_{n}\in[0,1];\quad I(t_{n}u_{n})=\max_{s\in[0,1]}I(su_{n}).$$
	Thus, given $R>0$, for $n$ large enough,
	$$I(t_{n}u_{n})\geq I\left(\frac{R}{\|u_{n}\|}u_{n}\right)=\frac{R^2}{2}-\frac{1}{2}\int_{\R^N}K(Rw_{n})F(Rw_{n})dz$$
	so,
	$$I(t_{n}u_{n})\geq \frac{R^2}{2}+o_{n}(1)$$
	from where it follows that,
	$$\liminf_{n\longrightarrow\infty}I(t_{n}u_{n})=\infty,$$
	as $I(0)=0$ and $I(u_{n})\longrightarrow c_{a}$, we have that $t_{n}\in(0,1)$. And so, $I'(t_{n}u_{n})(u_{n})=0$ and $I'(t_{n}u_{n})(t_{n}u_{n})=0$. Thus,
	$$4I(t_{n}u_{n})=4I(t_{n}u_{n})-I'(t_{n}u_{n})(t_{n}u_{n})=t_{n}^{2}\|u_{n}\|^{2}+\int_{\R^N}K(t_{n}u_{n})[f(t_{n}u_{n})(t_{n}u_{n})-2F(t_{n}u_{n})]dz$$
	consequently,
	$$4I(t_{n}u_{n})\leq\|u_{n}\|^{2}+\int_{\R^N}K(u_{n})[f(u_{n})(u_{n})-2F(u_{n})]dz=4I(u_{n})-I'(u_{n})(u_{n})=4I(u_{n})+o_{n}(1),$$
	a contradiction. Therefore, $(u_{n})$ is bounded.
	
\end{dem}

Note that, we have actually proved that, the Cerami sequence $(u_{n})$ is bounded in $E$, and so,  $u_{n}\rightharpoonup u$ in $E$. From where it follows that,
$$I'(u_{n})(u_{n})=o_{n}(1)\quad\text{and}\quad I'(u_{n})(u)=o_{n}(1)$$
so,
$$\|u_{n}\|^{2}=\int_{\R^N}K(u_{n})f(u_{n})u_{n}dz+o_{n}(1)$$
and
$$\|u\|^{2}=\int_{\R^N}K(u_{n})f(u_{n})udz+o_{n}(1).$$
Therefore,
$$\|u_{n}\|\longrightarrow\|u\|,$$
thence,
$$u_{n}\longrightarrow u\quad\text{em}\quad E.$$

Applying arguments above, easily, we have the \textbf{Theorem 1}.

\section{Periodic-Coercive Case}

Everything that was studied until the \textbf{Lemma 6} is hold for this new potential. Now, we prove a version of \textbf{Lemma 8} for this new potential.

\begin{lem}
	Let $(u_{n})\subset E_{a}$ such that $I(u_{n})\longrightarrow c_{a}$ and $I'(u_{n})(u_{n})\longrightarrow0$, then, $(u_{n})$ is bounded in $E_{a}$.
\end{lem}

\begin{dem}
	Firstly, note that, without loss of generality,  $u_{n}\geq0$ in $\R^N$, for all $n\in\N$. We claim that, $(u_{n})$ is bounded in $E_{a}$. Indeed, otherwise, $\|u_{n}\|\longrightarrow\infty$. Setting, $w_{n}=u_{n}/\|u_{n}\|$, we have
	$$w_{n}\rightharpoonup w\quad\text{in}\quad E_{a},$$
	with $w\leq0$. By \textbf{Lion's Lemma}, one of two situations occurs:
	
\begin{itemize}
	\item[$(i)$] For all $q\in(2,2^*)$,
	$$\lim_{n\longrightarrow\infty}\int_{\R^N}w_n^qdz=0.$$
	\item[$(ii)$] There exists $R,\eta>0$ and $(x_{n})\subset\Z^L$ such that,
	$$\liminf_{n\longrightarrow\infty}\int_{B_{R}(x_{n},0)}w^2_ndz\geq\eta>0.$$
\end{itemize}
	
In truth, the Lion's Lemma, ensures that $(i)$ or that for every $R>0$, there exists $\delta>0$ and $(z_n)\subset\R^N$, $z_n=(x_n,y_n)$, such that
$$\int_{B_R(x_n,y_n)}|w_n|^2dz\geq\delta>0,\quad\text{for all}\quad n\in\N.$$
But, note that,

\begin{Af}
	$(y_n)$ is bounded in $\R^M$.
\end{Af}
The proof of such a claim is similar to what we did at the beginning of \textbf{section 3}.

By the claim above, there exists $\lambda>0$ such that, $|y_n|\leq\lambda$, for all $n\in\N$. Then, clearly,
$$B_R(x_n,y_n)\subset B_{\tilde{R}}(x_n,0),\quad\text{for all}\quad n\in\N$$
where $\tilde{R}=R+\lambda$. Thus,
$$0<\delta\leq\int_{B_R(x_n,y_n)}|w_n|^2dz\leq\int_{ B_{\tilde{R}}(x_n,0)}|w_n|^2dz,$$
to facilitate, replace $\tilde{R}$ by $R$. Now, for all $n\in\N$, consider $\overline{x}_n\in\Z^L$ such that, $|\overline{x}_n-x_n|\leq \sqrt{L}$. Then, for $\overline{R}=R+\sqrt{L}$, we have
$$B_{R}(x_n,0)\subset B_{\overline{R}}(\overline{x}_n,0),$$
hence,
$$0<\delta\leq\int_{B_R(x_n,y_n)}|w_n|^2dz\leq\int_{B_{\overline{R}}(\overline{x}_n,0)}|u_n|^2dz,$$

It is then guaranteed that, indeed, $(i)$ or $(ii)$ is hold. If $(i)$ occurs, the result follows as the \textbf{Lemma 5}. Now, if $(ii)$ occurs, setting $\tilde{w}_n(z)=w_n(z+z_n)$, where $z_n=(x_n,0)$, next, $\|\tilde{w}_n\|=\|w_n\|$, so, $(\tilde{w}_n)$ is bounded in $E_a$. Thus, $\tilde{w}_n\rightharpoonup\tilde{w}$ in $E_a$ and then,
$$\liminf_{n\longrightarrow\infty}\int_{B_R(0)}\tilde{w}^2_ndz=\liminf_{n\longrightarrow\infty}\int_{B_R(z_n)}w^2_ndz\geq\eta>0$$
and so,
$$\int_{B_R(0)}\tilde{w}^2dz\geq\eta>0$$
and $\tilde{w}\neq0$. Now, as
$$I(u_{n})=c_{a}+o_{n}(1)\Longrightarrow\int_{\R^N}K(u_{n})\frac{F(u_{n})}{\|u_{n}\|^2}dz=1+o_n(1)$$
then,
$$\int_{\R^N}K(u_{n})\frac{F(u_{n})}{u_{n}^2}w_n^2dz\leq 1+o_n(1)\Longrightarrow\int_{B_{R}(z_n)}K(u_{n})\frac{F(u_{n})}{u_{n}^2}w_n^2dz\leq 1+o_n(1).$$
By changing variables and assumption of $f$, defining
$$G_n=\{z\in\R^N;\quad u_n(z+z_n)\geq\delta\},$$
we have for $M>0$, that
$$M\int_{B_{R}(0)\cap G_n}K(u_{n})(z+z_n)\tilde{w}_n^2dz\leq 1+o_n(1).$$
On the other hand, by $F(s)/s^2$ is increasing for $s>0$, we have
$$K(u_n)(z+z_n)\geq\|u_n\|^2\int_{\R^N}\frac{F(w_n(\overline{z}))}{|\overline{z}-(z+z_n)|^\gamma} d\overline{z}$$
thus,
$$K(u_n)(z+z_n)\geq\|u_n\|^2K(w_n)(z+z_n).$$
Therefore, by changing variables, $K(w_n)(z+z_n)=K(\tilde{w}_n)(z)$, thus,
$$1+o_n(1)\geq M\|u_n\|^2\int_{B_{R}(0)\cap G_n}K(\tilde{w}_n)(z)\tilde{w}_n^2dz=M\|u_n\|^2\int_{B_{R}(0)}\chi_{G_n}(z)K(\tilde{w}_n)(z)\tilde{w}_n^2dz.$$
By the properties of $F$,
$$K(\tilde{w}_n)(z)\geq\frac{F(\delta)}{\|u_{n}\|^2}\int_{\R^N}\chi_{B_{R}(0)\cap G_n}(\overline{z})|\overline{z}-z|^{-\gamma}d\overline{z},$$
thus,
	$$1+o_{n}(1)\geq M.F(\delta)\int_{B_R(0)}\chi_{G_n}(z)\tilde{w}_n^2(z)X_{n}(z)dz$$
	where,
	$$X_{n}(z)=\int_{\R^N}\chi_{B_{R}(0)\cap G_n}(\overline{z})|\overline{z}-z|^{-\gamma}d\overline{z}.$$
	For $z\in B_{R}(0)$, a. e., $z\in G_n$ for $n$ large and $\liminf_{n\longrightarrow\infty}X_{n}(z)>0$, we got that,
	$$1\geq M_{\delta}\int_{B_{R}(0)}\tilde{w}^2dz\geq M_{\delta}\eta>0.$$
	A contradiction. Therefore, $(u_{n})$ is bounded.	

\end{dem}

We got that, $(u_n)\subset E_a$ is bounded, so, $u_n\rightharpoonup u$ in $E_a$. Clearly, $I'(u)=0$, then, $u$ is solution of problem, not necessarily nontrivial. Thus, we need,

\begin{lem}
	There exist $(x_n)\subset\Z^L$, such that $(w_n)\subset E_a$, given by $w_n(x,y)=u_n(x+x_n,y)$, weakly converges to a function $w\in E_a$ with $w\neq0$ and $I'(w)=0$. In other words, $w$ is a solution nontrivial.
\end{lem}

\begin{dem}
	Firstly, we have the claim:
	\begin{Af}
		For all $R>0$, there exists $\delta>0$ and $(z_n)\subset\R^N$ such that
		$$\int_{B_R(z_n)}|u_n|^2dz\geq\delta>0,\quad\mbox{for all}\quad n\in\N.$$
	\end{Af}
	Indeed, otherwise, by \textbf{Lion's Lemma}, $u_n\longrightarrow0$ in $L^q(\R^N)$, for all $q\in(2,2^*)$. Thus, $u_n\longrightarrow0$ in $E_a$. A contradiction. Consequently, similarly to what we did in \textbf{Lemma 10}, there exist $R,\delta>0$ and $(x_{n})\subset\Z^L$ such that,
	$$\int_{B_{R}(x_{n},0)}u^2_ndz\geq\delta>0.$$
	Setting, $w_n(x,y)=u_n(x+x_n,y)$, we have $\|w_n\|=\|u_n\|$, and so, as $(u_n)$ is bounded, $(w_n)$ is bounded too. Hence, $w_n\rightharpoonup w$ in $E_a$, consequently, $w_n\longrightarrow w$ in $L^s_{loc}(\R^N)$, for all $s\in[2,2^*]$. Thus, $w_n\longrightarrow w$ in $L^2(B_R(0))$, but,
	$$0<\delta\leq\int_{B_{R}(x_{n},0)}|u_n|^2dz=\int_{B_{R}(0)}|w_n|^2dz$$
	and so, $w\neq0$.
	
	\begin{Af}
		$w$ is critical point of $I$.
	\end{Af}
	Indeed, note that
	$$I(w_n)=\frac{1}{2}A_n-\frac{1}{2}B_{n}$$
	where,
	$$A_n=\int_{\R^N}|\nabla w_n|^2+a(z)|w_n|^2dz\quad\text{and}\quad B_n=\int_{\R^N}K(w_n)F(w_n)dz.$$	
	By periodicity,
	$$A_n=\int_{\R^N}|\nabla u_n|^2+a(z)|u_n|^2dz$$
	and, by changing variables, we have
	$$B_n=\int_{\R^N}K(u_n)(z)F(u_n)(z)dz.$$
	Hence,
	$$I(w_n)=I(u_n)\longrightarrow c_a.$$
	On the other hand, for $\varphi\in E_a$ with $\|\varphi\|\leq1$, proceeding analogously, we have
	$$|I'(w_n)\varphi|=|I'(u_n)\varphi(z-z_n)|\leq\|I'(u_n)\|.\|\varphi\|\leq\|I'(u_n)\|$$
	where consider $z_n=(x_n,0)$ and , then, $\|I'(w_n)\|\longrightarrow0$ when $n\longrightarrow\infty$.
	Therefore, $w$ is a nontrivial critical point of problem.	
	
\end{dem}

By assumptions on $f$, we have $w\geq0$, thus, by the weak maximum principle, in \cite{Gilbarg-Trudinger}, we have  $w>0$.

We find the solution we wanted, but, we are not sure of being a ground state solution. But, note that, $w_n\rightharpoonup w$ in $E_a$, and so
$$w_n(z)\longrightarrow w(z)\quad\text{and}\quad K(w_n)f(w_n)w_n(z)\longrightarrow K(w)f(w)w(z)\quad\text{a.e. in $\R^N$.}$$
Moreover,
$$\|w\|^2=\int_{\R^N}K(w)f(w)wdz,$$
thus, by Fatou's Lemma,
$$I(w)=I(w)-\frac{1}{2}I'(w)(w)\leq \varliminf_{n\longrightarrow\infty}\left(I(w_n)-\frac{1}{2}I'(w_n)w_n\right)=c_a,$$
therefore, $I(w)\leq c_a$.

To be $w$ ground state solution, we have the theorem below:

\begin{teo}
	Let $N=\{u\in E_a;\quad I'(u)(u)=0\}\setminus\{0\}$, the Nehari manifold of $I$. Then,
	$$c_a\leq\inf_{N}I.$$
\end{teo}

\begin{dem}
	Let $u\in N$. For $t\geq0$, we got
	$$I(tu)=\frac{t^2}{2}\|u\|^2-\frac{1}{2}\int_{\R^N}K(tu)F(tu)dz$$	
	and
	$$0=I'(u)(u)=\|u\|^2-\int_{\R^N}K(u)f(u)udz.$$
	Moreover, using properties of $f$, we have that
    $$\frac{d}{dt}[I(tu)]>0,\quad\text{if $t<1$ and}\quad\frac{d}{dt}[I(tu)]<0,\quad\text{if $t>1$}.$$	
	Therefore,
	$$\max_{t\geq0}I(tu)=I(u)$$
	and so,
	$$\inf_{\alpha\in\Gamma}\max_{t\in[0,1]}I(\alpha(t))\leq I(u),\quad \forall u\in N.$$
	Consequently,
	$$c_a\leq\inf_{N}I.$$
\end{dem}

This proves the existence of a solution in the mountain pass level, i.e., thus proving the \textbf{theorem 2}.
 
\section{Asymptotically Periodic-Coercive Case}

Note that, associated with potential $a_{p}$, we have the problem:

\[
(Q)\left\{
\begin{array}{lcl}
-\Delta u+a_p(z)u =K(u)f(u),\quad \R^{N} \\
u\in H^{1}(\R^{N}),\quad u>0\quad\text{in $\R^{N}$} 
\end{array}
\right.
\]

As in the problem $(P)$, we must find solution in a space of type $E_a$, i.e., in the  $E_{a_p}$, analogous to $E_a$. Now, note that, $E_a=E_{a_p}$ and, $\|.\|_{E_a}$ and $|.\|_{E_{a_p}}$ are equivalents.

Consequently, we have that, the energy functional of problems $(P)$ and $(Q)$ are $I,J:E\longrightarrow\R$, defined respectively by:
$$I(u)=\frac{1}{2}\int_{\R^N}|\nabla u|^2+a(z)|u|^2dz-\frac{1}{2}\int_{\R^N}K(u)F(u)dz$$
and
$$J(u)=\frac{1}{2}\int_{\R^N}|\nabla u|^2+a_p(z)|u|^2dz-\frac{1}{2}\int_{\R^N}K(u)F(u)dz.$$

These energy functionals check up the mountain pass geometry. For $I$, exist $(u_n)\subset E$ such that,
$$I(u_n)\longrightarrow c_a\quad\text{and}\quad I'(u_n)(u_n)\longrightarrow0.$$

Similarly to what we have already done, $(u_n)$ is bounded in $E$, with any of the norms $\|.\|_{E_a}$ and $|.\|_{E_{a_p}}$. And so, $u_n\rightharpoonup u$ in $E$, $u_n(z)\longrightarrow u(z)$ a.e. in $\R^N$ and $u_n\longrightarrow u$ in $L^s_{loc}(\R^N)$, for $s\in[1,2^*)$. Moreover, $I'(u)=0$. We need to ensure that $u\neq0$.

\begin{Af}
	$u\neq0$.
\end{Af}

By the previous case, there exists $w\in E$ positive, such that
$$J(w)=c_{a_p}\quad\text{and}\quad J'(w)=0$$
where $c_{a_p}$ is the level of Mountain Pass for $J$.

\begin{lem}
	The levels $c_a$ and $c_{a_p}$ satisfies $c_a<c_{a_p}$.
\end{lem}

\begin{dem}
	Clearly, by definition, we have that $c_a\leq c_{a_p}$. For $w$, consider the path
	$$\delta_{w}(t)=tsw,\quad\text{for all}\quad t\in[0,1]\quad(\delta_{w}\in\Gamma)$$
	where $s$ is fixed, such that
	$$I(sw)<0.$$
	So,
	$$c_a\leq\max_{t\in[0,1]}I(\delta_{w}(t))=\max_{t\in[0,1]}I(tsw)\leq\max_{r\geq0}I(rw)=I(r_{0}w)<J(r_{0}w)\leq\max_{r\geq0}J(rw)=J(w),$$
	thus, 
	$$c_a<J(w)=c_{a_p}.$$	
\end{dem}

Now, note that, if $u=0$, we have:

\begin{lem}
	$$|I(u_n)-J(u_n)|\longrightarrow0\quad\text{and}\quad \|I'(u_n)-J'(u_n)\|\longrightarrow0.$$
\end{lem}

It is very important to highlight, as $(u_n)$ is bounded in $E$, that there exists $M>0$ such that $\|u_n\|\leq M$, for all $n\in\N.$ Thus, as for all $\varepsilon>0$, there exists $R>0$, such that
$$a(x,y)\geq\frac{1}{\varepsilon},\quad |y|\geq R,$$
we have,
$$\int_{|y|\geq R}|u_n|^2dz\leq\int_{|y|\geq R}\varepsilon a(z)|u_n|^2dz\leq\varepsilon\int_{|y|\geq R}a(z)|u_n|^2dz\leq \varepsilon M,$$
for interpolation,
$$\int_{|y|\geq R}|u_n|^qdz<C\varepsilon,\quad\mbox{for all}\quad q\in[2,2*)$$
and consequently,
$$\left|\int_{|y|\geq R}K(u_n)f(u_n)u_ndz\right|<C\varepsilon$$
where the above properties are valid for all $n\in\N$, and the constant $C>0$. 

For $R>0$ large enough, define $\varphi_{R}\in C^\infty(\R^M)$, such that $|\nabla\varphi_R|\leq2/R$, $\varphi_{R}(y)=0$ if $|y|\leq R/2$, $\varphi_{R}(y)=1$ if $|y|\geq R$ and $0\leq\varphi_{R}(y)\leq1$, for all $y\in\R^M$. We define also, $v_n(x,y)=u_n(x,y)\varphi_{R}(y)$, so, $\nabla v_n=(\nabla u_n)\varphi_{R}+u_n(\nabla\varphi_{R})$, clearly, $(v_n)$ is bounded in $E$ and, consequently, we have that
$$o_{n}(1)=I'(u_n)(v_n)=\int_{\R^N}|\nabla u_n|^2\varphi_R+a(z)u_n^2\varphi_{R}dz+\int_{\R^N}u_n\nabla\varphi_R\nabla u_ndz-\int_{\R^N}K(u_n)f(u_n)u_n\varphi_Rdz.$$
Noting that, for $\varepsilon>0$ small enough and $R>0$ large enough, we have
$$\left|\int_{\R^N}u_n\nabla\varphi_R\nabla u_ndz\right|<\varepsilon$$
as well as,
$$\left|\int_{\R^N}K(u_n)f(u_n)u_n\varphi_Rdz\right|<\varepsilon.$$
Thus, for $\varepsilon>0$ small enough and $R>0$ large enough
$$\int_{\R^N}|\nabla u_n|^2\varphi_R+a(z)u_n^2\varphi_{R}dz<\varepsilon$$
so,
$$\int_{|y|\geq R}a(z)|u_n|^2dz<\varepsilon.$$
Similarly, you can check that
$$\int_{|y|\geq R}a_{p}(z)|u_n|^2dz<\varepsilon.$$

Now, we can prove the \textbf{lemma 13}:

\begin{dem}
	Indeed,
	\begin{eqnarray*}
		I(u_n)-J(u_n)&=&\int_{\R^L\times \overline{B}_{R}^c(0)}[a(z)-a_{p}(z)]|u_n|^2dz+\int_{\overline{B}_{R}(0)\times \overline{B}_{R}(0)}[a(z)-a_{p}(z)]|u_n|^2dz\\
		&&+\int_{\overline{B}_{R}^c(0)\times \overline{B}_{R}(0)}[a(z)-a_{p}(z)]|u_n|^2dz
	\end{eqnarray*}
		thus, by asymptotic properties,
		$$\left|\int_{\R^L\times \overline{B}_{R}^c(0)}[a(z)-a_{p}(z)]|u_n|^2dz\right|<\frac{\varepsilon}{3}$$
		as, $u_n\longrightarrow0$ in $L^2(\overline{B}_{R}(0)\times \overline{B}_{R}(0))$ and the limitation of the potential,
		$$\left|\int_{\overline{B}_{R}(0)\times \overline{B}_{R}(0)}[a(z)-a_{p}(z)]|u_n|^2dz\right|<\frac{\varepsilon}{3}$$
		and, 
		$$\left|\int_{\overline{B}_{R}^c(0)\times \overline{B}_{R}(0)}[a(z)-a_{p}(z)]|u_n|^2dz\right|<\frac{\varepsilon}{3}$$
		we got that,
		$$|I(u_n)-J(u_n)|\longrightarrow0.$$
		Similarly,
		$$\|I'(u_n)-J'(u_n)\|\longrightarrow0.$$

\end{dem}

In this way, we ensure that,
$$J(u_n)\longrightarrow c_a\quad\text{and}\quad J'(u_n)(u_n)\longrightarrow0.$$
By \textbf{Lema 7}: exist $(\overline{x}_n)\subset\Z^L$, such that $(v_n)\subset E$, given by $v_n(x,y)=u_n(x+\overline{x}_n,y)$, converges weakly to $v\in E$ with $v\neq0$ and $J'(v)=0$. Thus, $w$ is solution nontrivial. And so, $v\in N_{J}=\{u\in E;\quad J'(u)(u)=0\}\setminus\{0\}$, in this way,
$$c_{a_p}\leq\inf_{u\in N_{J}}J(u)\leq J(v).$$
On the other hand,
\begin{eqnarray*}
	J(v)&=&J(v)-\frac{1}{2}J'(v)(v)\\
		&\leq&\varliminf_{n\longrightarrow\infty}\int_{\R^N}K(v_n)\frac{f(v_n)v_n}{2}-K(v_n)F(v_n)dz\\
		&=&\varliminf_{n\longrightarrow\infty}\left(I(u_n)-\frac{1}{2}I'(u_n)(u_n)\right)
\end{eqnarray*}
therefore,
$$J(v)\leq\varliminf_{n\longrightarrow\infty}I(u_n)=c_a,$$
so, $c_{a_p}\leq c_a$. A contradiction. Therefore, $u\neq0$.

By applying the arguments above, we have \textbf{theorem 3}.

\end{document}